\documentclass[11pt]{article}
\usepackage{amsmath, amssymb, amsthm}

\theoremstyle{plain}

\numberwithin{equation}{section}

\title{On higher-power moments of $\Delta(x)$\hspace{2mm}(III) }
\date{}

\begin{document}

\maketitle

\vspace{-1.5cm}

\begin{center}

Wenguang  Zhai

\medskip

School of Mathematical Sciences \\
Shandong Normal University \\
Jinan 250014, Shandong \\
 P. R. China\\
zhaiwg@hotmail.com
\end{center}

\begin{center}{Acta Arith.
{\bf 118} (2005), 263--281}\end{center}

\begin{abstract}
 Let $\Delta(x)$ be the error term of the Dirichlet
divisor problem.   An asymptotic formula with the error term
$O(T^{53/28+\varepsilon})$ is established for the integral
$\int_1^T\Delta^4(x)dx.$
  Similar results are also established for some
other well-known error terms in the analytic number theory .
\end{abstract}
   \footnote[0]{This work is
supported by National Natural Science Foundation of China(Grant No.
10301018).}

\section{\bf  Introduction and  main results }

Let $d(n)$ denote the Dirichlet divisor function and $\Delta(x)$
denote the error term of the sum  $ \sum_{n\leq x}d(n)$ for a large
real variable $x.$ Dirichlet first proved that
 $\Delta(x)=O(x^{1/2}).$
 The exponent $1/2$ was improved by many authors.
The latest result reads
\begin{equation}
\Delta(x)\ll x^{131/416}(\log x)^{26947/8320 },
\end{equation}
 proved by  Huxley[3]. It
is conjectured that
\begin{equation}
\Delta(x)=O(x^{1/4+\varepsilon}),
\end{equation}
 which is supported by the
classical mean-square  result
\begin{equation}
\int_1^T\Delta^2(x)dx=\frac{(\zeta(3/2))^4}{6\pi^2 \zeta(3)}T^{3/2}
+O(T\log^5 T)
\end{equation}
proved by Tong[10].

Tsang[11] studied the third- and fourth-power moments of
$\Delta(x).$ He proved that the asymptotic formulas
\begin{equation}
\int_2^T\Delta^3(x)dx=\frac{3c_1}{28\pi^3}T^{7/4}+O(T^{7/4-\delta_1+\varepsilon}
)
\end{equation}
and
\begin{equation}
\int_2^T\Delta^4(x)dx=\frac{3c_2}{64\pi^4}T^2+O(T^{2-\delta_2+\varepsilon})
\end{equation}
hold, where $\delta_1=1/14,$ $\delta_2=1/23,$
\begin{eqnarray*}
&&c_1:=\sum_{\alpha,\beta,h\in {\Bbb N}}
(\alpha\beta(\alpha+\beta))^{-3/2}h^{-9/4}|\mu(h)|d(\alpha^2
h)d(\beta^2 h)
d((\alpha+\beta)^2 h),\\
&&c_2:=\sum_{\stackrel{n,m,k,l\in {\Bbb N}}{\sqrt n+\sqrt m=\sqrt
k+\sqrt l}} (nmkl)^{-3/4}d(n)d(m)d(k)d(l).
\end{eqnarray*}

Recently in [12] the author  proved that (1.4) holds  for
$\delta_1=1/4.$   Ivi\'c and Sargos[7] proved that (1.4) holds for
$\delta_1=7/20.$ The author got this exponent independently. However
, Professor Ivi\'c kindly informed the author that the exponent
$\delta_1=7/20$ had already been obtained by Professor Tsang several
years ago but he had never published this result.

 Following Tsang's approach , in [12] the author proved that
 (1.5) holds for  $\delta_2=2/41.$ This approach used the method of exponential sums.
Especially if the exponent pair conjecture is true, namely, if
$(\varepsilon,1/2+\varepsilon)$ is an exponent pair, then (1.5)
holds for $\delta_2=1/14.$ However, in [7]
 Ivi\'c and Sargos ingeniously proved a substantially better result .
 They proved that (1.5) holds for
   $ \delta_2=1/12.$

In this paper,  combining the method of [7] and a recent deep result
of Robert and Sargos[9], we shall prove the following

{\bf Theorem  1. } We have
\begin{equation} \int_2^T\Delta^4(x)dx=\frac{3c_2}{64\pi^4}
T^2+O(T^{53/28+\varepsilon}).
\end{equation}

The theorem is also true for other error terms . Let $P(x)$ denotes
the error term of the Gauss circle problem , which is an error term
similar to $\Delta(x).$ Let $a(n)$ be the Fourier coefficients of a
holomorphic cusp form of weight $\kappa=2n\geq 12$ for the full
modular group and define
$$A(x):=\sideset{}{^{\prime}}\sum_{n\leq x}a(n),\hspace{3mm}x\geq 2.$$
We then have the following two corollaries, which improve the
previous results([2], [11], [12]).

 {\bf Corollary 1.} We
have
\begin{equation}
\int_2^TP^4(x)dx=CT^{2}+O(T^{53/28+\varepsilon}).
\end{equation}

{\bf Corollary 2. } We have
\begin{equation}
\int_1^T A^4(x)dx=B_\kappa T^{2\kappa}+O(T^{2\kappa-
3/28+\varepsilon}).
\end{equation}

\bigskip

Now let's consider $E(t),$ defined by
\begin{equation}
E(t):=\int_0^t|\zeta(\frac{1}{2}+iu)|^2du-t\log(t/2\pi)-(2\gamma-1)t,\hspace{2mm
}t\geq 2.
\end{equation}

Tsang[11] also studied the fourth-power moment of $E(t)$ by using
Atkinson's formula[1] and proved that
\begin{equation}
\int_2^TE^4(t)dt=\frac{3}{8\pi}c_2T^2 +O(T^{2-\delta_3+\varepsilon})
\end{equation}
with some unspecified constant $\delta_3>0.$

  Ivi\'c[4] used a different way to study the higher power moments of $E(t)$.
   The following is his approach. Let
\begin{equation}
\Delta^{*}(x):= \frac{1}{2}\sum_{n\leq 4x}(-1)^nd(n)-x(\log
x+2\gamma-1),\hspace {2mm} x\geq 1.
\end{equation}
Then for $1\ll N\ll x,$ we have[6]
\begin{equation}
\Delta^{*}(x)=\frac{1}{\pi\sqrt 2}\sum_{n\leq
N}(-1)^nd(n)n^{-3/4}x^{1/4} \cos(4\pi\sqrt{nx}-\pi/4)
+O(x^{1/2+\varepsilon}N^{-1/2})
\end{equation}

Jutila[8] proved that
\begin{equation}
\int_0^T(E(t)-2\pi\Delta^{*}(\frac{t}{2\pi}))^2dt\ll T^{4/3}\log^3
T,
\end{equation}
which means that $E(t)$ is well approximated by
$2\pi\Delta^{*}(\frac{t}{2\pi})$ at least in the mean square sense.
From (1.13) Ivi\'c[4] deduced  that
\begin{equation}
\int_0^TE^4(t)dt=(2\pi)^5\int_0^{\frac{T}{2\pi}}(\Delta^{*}(t))^4dt
+O(T^{23/12}\log^{3/2} T).
\end{equation}
Thus the fourth power moment of $E(t)$ was transformed into the
fourth power moment of $\Delta^{*}(t),$  which can be dealt with in
the same way as the fourth power moment of $\Delta(x).$ By Tsang's
result[11] , Ivi\'c deduced from (1.14) that (1.10) holds for
$\delta_3=1/23.$ In [7], Ivi\'c and Sargos  proved that  one can
take $\delta_3=1/12.$

It is easy to see that $1/12$ is the limit of this approach  since
it is the limit of Jutila's result (1.13). In this paper, we shall
use a different way to prove the following

{\bf Theorem 2.} We have
\begin{equation}
\int_2^TE^4(t)dt=\frac{3}{8\pi}c_2T^2 +O(T^{53/28+\varepsilon}).
\end{equation}

{\bf Remark.} The proof of Theorem 2 doesn't use (1.13) and it is
actually a generalization of the approach used in the author [13].
In [14] the author used a similar method to study the third power
moment of $E(t).$

 {\bf Acknowledgement.} I am   very grateful
to the referee for many valuable suggestions . Especially he gave a
direct proof of  Lemma 3, which is much easier than my previous
approach. I wish to
 thank Professor Ivi\'c and Professor Robert, who kindly send me
the papers [7] and [9], respectively.

{\bf Notations.} Throughout
 this paper,
$[x]$ denotes the integer part of $x,$ $\Vert x\Vert$ denotes the
distance from $x$ to the integer nearest to $x,$ $n\sim N$ means
$N<n\leq 2N,$  $n\asymp N$ means $C_1N<n\leq C_2N$ for positive
constants $C_1<C_2.$ $\varepsilon$ always denotes a small positive
constant which may be different at different places. We shall use
the estimate  $d(n)\ll n^\varepsilon$ freely.

\section{\bf The spacing problem of the square roots}

In the proofs of Theorem 1 and 2 , the sums and differences of
square roots will appear in the exponential. Thus we should study
the spacing problem of the square roots.

We need the following Lemmas. Lemma 1 is a special case of  a new
result proved in Robert and Sargos[9], which also plays an important
role in this paper. Lemma 2 is Lemma 3 of Tsang [11].
 Lemma 3 provides an upper bound
of the number of solutions of the inequality
\begin{equation}
|n_1^{1/2}+n_2^{1/2}\pm
n_3^{1/2}-n_4^{1/2}|<\Delta,\hspace{5mm}n_j\sim N_j(j=1,2,3,4),
\end{equation}
where $N_j\geq 2(j=1,2,3,4)$ are real numbers . Lemma 4 is
essentially  Lemma 3 of Ivi\'c and Sargos[7], but we added the case
$\alpha\ll 1.$ Lemma 5 is essentially Lemma 5 of [7], but the term
$K\min(M,M^{\prime},L)$ therein is superfluous since we add the
condition $|\sqrt n+\sqrt m-\sqrt k-\sqrt l|>0$ in Lemma 5, and so
we give a new proof here. Lemma 6 is Lemma 6 of [7].

 {\bf Lemma 1.} Suppose $N\geq 2,$ $\Delta>0.$
Let ${\cal A}(N;\Delta)$ denote the number of solutions of the
inequality $$|n_1^{1/2}+n_2^{1/2}-
n_3^{1/2}-n_4^{1/2}|<\Delta,\hspace{5mm}n_j\sim N(j=1,2,3,4),$$ then
$${\cal A}(N;\Delta)\ll (\Delta N^{7/2}+N^{2})N^\varepsilon.$$

{\bf Lemma 2.} If $n,m,k,l\in {\Bbb N}$  such that $\sqrt n+\sqrt
m\pm\sqrt k-\sqrt l\not= 0,$
 then respectively,
$$|\sqrt n+\sqrt m\pm\sqrt k-\sqrt l
|\gg \max(n,m,k,l)^{-7/2}.$$

{\bf Lemma 3. } Suppose $N_j\geq 2(j=1,2,3,4),$ $\Delta>0.$ Let
${\cal A}_{\pm}(N_1,N_2,N_3,N_4;\Delta)$ denote the number of
solutions of the inequality (2.1),  then we have
$${\cal A}_{\pm}(N_1,N_2,N_3,N_4;\Delta)\ll
\prod_{j=1}^4(\Delta^{1/4} N_j^{7/8}+N_j^{1/2})N_j^{\varepsilon}.$$

\begin{proof} We shall use a combinatorial argument to prove this
Lemma. Let $\{a_i\}$ and $\{b_i\}$  be two finite sequences of real
numbers. Let $\Delta>0$ . Suppose $u_0$ and $J$(a positive integer)
are chosen so that $\{a_i\}\subset (u_0,u_0+J\Delta],$
 $\{b_i\}\subset (u_0,u_0+J\Delta].$ Divide this interval into the
 abutting subintervals $I_j:=(u_0+j\Delta,u_0+(j+1)\Delta]$
for $j=0,1,\cdots, J-1$ and then let
$$N_j(A):=\#\{i:a_i\in I_j\}, N_j(B):=\#\{i:b_i\in I_j\}. $$
If $|a_r-b_s|\leq \Delta$, then either both $a_r$ and $b_s$ lie in
the same subinterval $I_j,$ or they lie in adjacent subintervals
$I_j$ and $I_{j+1}.$ Hence
\begin{eqnarray*}
&&\#\{(r,s): |a_r-b_s|\leq \Delta\}\\
&&\leq
\sum_{j}N_{j}(A)N_{j}(B)+\sum_{j}N_{j}(A)N_{j+1}(B)+\sum_{j}N_{j+1}(A)N_{j}(B)\\
&&\leq 3(\sum_{j}N_{j}(A)^2)^{1/2}(\sum_{j}N_{j}(B)^2)^{1/2}
\end{eqnarray*}
by Cauchy-Schwarz's inequality. On the other hand , we have
\begin{eqnarray*}
&&\sum_{j}N_{j}(A)^2=\sum_{j}\#\{(r,r^{\prime}):a_r,a_{r^{\prime}}\in
I_j\}\\
&&\leq \#\{(r,r^{\prime}):|a_r-a_{r^{\prime}}|\leq \Delta\},
\end{eqnarray*}
and similarly for $\sum_{j}N_{j}(B)^2.$ Thus we get the bound
\begin{equation}
\#\{(r,s): |a_r-b_s|\leq \Delta\}
\end{equation}
$$\leq 3 (\#\{(r,r^{\prime}):|a_r-a_{r^{\prime}}|\leq \Delta\})^{1/2}
(\#\{(s,s^{\prime}):|b_s-b_{s^{\prime}}|\leq \Delta\} )^{1/2}.$$

Suppose $n_j,n_j^{\prime}\sim N_j(j=1,2,3,4).$ Applying (2.2) to the
sequences $A=\{\sqrt{n_1}+\sqrt {n_2}\}$ and $B=\{\sqrt{n_3}+\sqrt
{n_4}\}$, we get
\begin{equation}
{\cal
A}_{-}(N_1,N_2,N_3,N_4)=\#\{(n_1,n_2,n_3,n_4):|n_1^{1/2}+n_2^{1/2}-n_3^{1/2}-n_4^{1/2}|\leq
\Delta\}
\end{equation}
\begin{eqnarray*}
&&\leq
3(\#\{(n_1,n_2,n_1^{\prime},n_2^{\prime}):|n_1^{1/2}+n_2^{1/2}-n_1^{\prime
1/2}-n_2^{\prime 1/2}|\leq \Delta\})^{1/2}\\
&&\hspace{4mm}\times
(\#\{(n_3,n_4,n_3^{\prime},n_4^{\prime}):|n_3^{1/2}+n_4^{1/2}-n_3^{\prime
1/2}-n_4^{\prime 1/2}|\leq \Delta\})^{1/2}.
\end{eqnarray*}
Using the previous bound to the sequences
$A_1=\{n_1^{1/2}-n_1^{\prime 1/2}\}, B_1=\{n_2^{1/2}-n_2^{\prime
1/2}\},$  and  $A_2=\{n_3^{1/2}-n_3^{\prime 1/2}\},
B_2=\{n_4^{1/2}-n_4^{\prime 1/2}\},$ respectively,  we get
\begin{equation}
{\cal A}_{-}(N_1,N_2,N_3,N_4)\leq 9\prod_{j=1}^4{\cal
A}_{-}(N_j,N_j,N_j,N_j)^{1/4},
\end{equation}
which combining Lemma 1 gives Lemma 3 for the case "-". The proof
for the case "+" is similar.

\end{proof}

{\bf Lemma 4.} Suppose $K\geq 10,$
 $\alpha, \beta\in {\Bbb R},$ $2K^{-1/2}\leq |\alpha|\ll K^{1/2}$ and $0<\delta<1/2.$ Then
  we have $$
\#\{k\sim K: \Vert \beta+\alpha\sqrt k\Vert <\delta\} \ll
K\delta+K^{1/2+\varepsilon}.$$

\begin{proof} Without loss of generality, suppose $\alpha>0.$
Let ${\cal N}=\#\{k\sim K: \Vert \beta+\alpha\sqrt k\Vert
<\delta\}.$  If $1\ll \alpha \ll K^{1/2},$ from Lemma 3 of Ivi\'c
and Sargos[7] we get
$${\cal N}\ll
K\delta+|\alpha|^{1/2}K^{1/4+\varepsilon}+K^{1/2+\varepsilon}\ll
K\delta+K^{1/2+\varepsilon}.$$

Now suppose $2K^{-1/2}\leq\alpha\ll 1.$ Since $\Vert t\Vert$ is a
periodic function with period $1,$ we suppose $0<\beta\leq 1.$ If
 $\Vert \beta+\alpha\sqrt k\Vert <\delta $, then
 there exists a unique $l\in  [\alpha\sqrt K, 2\alpha\sqrt K+2]$
such that
$$ (l-\beta-\delta)^2/\alpha^2<k\leq (l-\beta+\delta)^2/\alpha^2,$$
which implies
\begin{eqnarray*}
{\cal N}&&\ll \sum_{l\sim \alpha\sqrt
K}\left([(l-\beta+\delta)^2/\alpha^2]-[(l-\beta-\delta)^2/\alpha^2]\right)\\
&&\ll \sum_{l\sim \alpha\sqrt
K}\left((l-\beta+\delta)^2/\alpha^2-(l-\beta-\delta)^2/\alpha^2+1\right)\\
&&\ll K\delta+K^{1/2}
\end{eqnarray*}
if we notice $\alpha\ll 1.$
 \end{proof}

 {\bf Lemma 5.} Suppose $1\leq N\leq M, 1\leq L\leq K,$ $N\leq L,$
$M\asymp K,$ $0<\Delta\ll K^{1/2}.$ Let ${\cal A}_1(N,M,K,L;\Delta)$
denote the number of solutions of the inequality $$ 0<|\sqrt n+\sqrt
m-\sqrt k-\sqrt l|<\Delta$$ with $n\sim N, m\sim M, k\sim K, l\sim
L.$ Then  we have
$${\cal A}_1(N,M,K,L;\Delta)\ll \Delta K^{1/2}NML
+NLK^{1/2+\varepsilon}.$$ Especially if  $\Delta K^{1/2}\gg 1,$ then
$${\cal A}_1(N,M,K,L;\Delta)\ll \Delta K^{1/2}NML.$$

\begin{proof}
If $(n,m,k,l)$ satisfies $|\sqrt n+\sqrt m-\sqrt k-\sqrt l|<\Delta$
, then we get
$$m=k+2k^{1/2}(\sqrt l-\sqrt n)+(\sqrt l-\sqrt n)^2+u$$
with $|u|\leq C\Delta K^{1/2}$ for some absolute constant $C>0.$
Hence ${\cal A}_1(N,M,K,L;\Delta)$ does not exceed the number of
solutions of the inequality
\begin{equation}
|2k^{1/2}(\sqrt l-\sqrt n)+(\sqrt l-\sqrt n)^2+k-m|<C\Delta K^{1/2}
\end{equation}
with $n\sim N, m\sim M, k\sim K, l\sim L.$

If $\Delta K^{1/2}\gg 1,$ then for fixed $(n,k,l),$ the number of
$m$ for which (2.5) holds is $\ll 1+\Delta K^{1/2}\ll \Delta
K^{1/2}$ if we notice $K\asymp M.$ Hence $${\cal
A}_1(N,M,K,L;\Delta)\ll \Delta K^{1/2}NML.$$

Now suppose $\Delta K^{1/2}\leq 1/4C.$ For fixed $(n,k,l),$ there is
at most one $m$ such that (2.14) holds. If such $m$ exists, then we
have
\begin{equation}
\Vert 2k^{1/2}(\sqrt l-\sqrt n)+(\sqrt l-\sqrt n)^2 \Vert < C\Delta
K^{1/2}.\end{equation}

We shall use  Lemma 4 to bound the number of solutions of (2.6) with
$\alpha=2(\sqrt l-\sqrt n),\beta=(\sqrt l-\sqrt n)^2.$ Let ${\cal
C}_1$ denote the number of solutions of (2.6) with $|\alpha|\geq
2K^{-1/2},$ and ${\cal C}_2$ the number of solutions with $|\alpha|<
2K^{-1/2},$ respectively. By Lemma 4 we get
$${\cal C}_1\ll \Delta K^{1/2}NML+NLK^{1/2+\varepsilon}$$
if we notice $M\asymp K.$ Now we estimate ${\cal C}_2.$ From
$|\alpha|< 2K^{-1/2},$ we get $N\asymp L.$ If $l=n, $ then from
(2.5) we get $k=m.$ This contradicts to $|\sqrt n+\sqrt m-\sqrt
k-\sqrt l|>0.$ Thus $l\not= n. $ From
$$2K^{-1/2}> |\sqrt l-\sqrt n|=\frac{|l-n|}{\sqrt l+\sqrt
n}\geq \frac{1}{\sqrt l+\sqrt n}\geq 1/2\sqrt{2L}$$ we get $L\gg K$
and thus $N\asymp M\asymp K\asymp L.$  So we have $${\cal C}_2\ll
\#\{(l,n): |\alpha|< 2K^{-1/2}\}\times \#\{k\}\ll K^2,$$ which can
be absorbed into the estimate of ${\cal C}_1.$ This completes the
proof of Lemma 5.
\end{proof}

{\bf Lemma 6.} Suppose $1\leq N\leq M\leq K\asymp L,$ $0<\Delta\ll
L^{1/2}.$ Let ${\cal A}_2(N,M,K,L;\Delta)$ denote the number of
solutions of the inequality $$ |\sqrt n+\sqrt m+\sqrt k-\sqrt
l|<\Delta$$ with $n\sim N, m\sim M, k\sim K, l\sim L.$ Then  we have
$${\cal A}_2(N,M,K,L;\Delta)\ll \Delta L^{1/2}NMK
+NMK^{1/2+\varepsilon}.$$ Especially if  $\Delta L^{1/2}\gg 1,$ then
$${\cal A}_2(N,M,K,L;\Delta)\ll \Delta L^{1/2}NMK.$$

\section{\bf Proof of  Theorem 1}

 Suppose $T\geq 10.$ It suffices to evaluate the integral
 $\int_T^{2T}\Delta^4(x)dx.$
Suppose $y=T^{3/4}.$ For any $T\leq x\leq 2T,$  by the truncated
Voronoi's formula we get
\begin{equation}
\Delta(x)=\frac{1}{\sqrt 2\pi}{\cal R
}+(x^{1/2+\varepsilon}y^{-1/2}),
\end{equation}
where
\begin{eqnarray*}
&&{\cal R}:={\cal R}(x)=x^{1/4}\sum_{n\leq y}\frac{d(n)}
{n^{3/4}}\cos(4\pi\sqrt{xn}-\frac{\pi}{4}).
\end{eqnarray*}
 We have
\begin{align}
\int_T^{2T}\Delta^4(x)dx
 &=\frac{1}{4\pi^4}\int_T^{2T}{\cal
 R}^4dx+O(T^{9/4+\varepsilon}y^{-1/2}+T^{3+\varepsilon}y^{-2})\\
 &=\frac{1}{4\pi^4}\int_T^{2T}{\cal
 R}^4dx+O(T^{15/8+\varepsilon}).\nonumber
\end{align}

Let
$$g=g(n,m,k,l):=
(nmkl)^{-\frac{3}{4}}d(n)d(m)d(k)d(l),
\mbox{for}\hspace{2mm}n,m,k,l\leq y,$$ and $g=0$ otherwise.

The equation (3.4) of Tsang[11] reads
\begin{equation}
{\cal R}_1^4=S_1(x)+S_2(x)+S_3(x)+S_4(x),
\end{equation}
where
\begin{eqnarray*}
&&S_1(x):=\frac{3}{8}
\sum_{\sqrt n+\sqrt m=\sqrt k+\sqrt l}g x,\\
&&S_2(x):=\frac{3}{8} \sum_{\sqrt n+\sqrt m\not=\sqrt k+\sqrt l}g
x\cos(4\pi(\sqrt n+\sqrt m-\sqrt k-\sqrt l)\sqrt x),\\
&&S_3(x):=\frac{1}{2} \sum g x
\sin(4\pi(\sqrt n+\sqrt m+\sqrt k-\sqrt l)\sqrt x),\\
&&S_4(x):=-\frac{1}{8} \sum g x \cos(4\pi(\sqrt n+\sqrt m+\sqrt
k+\sqrt l)\sqrt x).
\end{eqnarray*}

From (3.7) of [11] we get
\begin{equation}
\int_T^{2T}S_1(x)dx=\frac{3c_2}{8}\int_T^{2T}xdx+
O(T^{2-3/16+\varepsilon}).
\end{equation}

From the first derivative test  we get
\begin{equation}
\int_T^{2T}S_4(x)dx\ll T^{3/2+\varepsilon}y^{1/2}\ll
T^{15/8+\varepsilon}.
\end{equation}

Now let us consider the contribution of $S_2(x).$ By the first
derivative test we get
\begin{equation}
\int_T^{2T}S_2(x)dx\ll \sum_{\stackrel{n,m,k,l\leq y} {\sqrt n+\sqrt
m\not=\sqrt k+\sqrt l}} g\min \left(T^2,\frac{T^{\frac{3}{2}}}
{|\sqrt n+\sqrt m-\sqrt k-\sqrt l|}\right)
\end{equation}
$$\ll T^{\varepsilon}G(N,M,K,L),$$
where
\begin{eqnarray*}
&&G(N,M,K,L)=\sum_1 g\min \left(T^2,\frac{T^{\frac{3}{2}}}
{|\sqrt n+\sqrt m-\sqrt k-\sqrt l|}\right),\\
&&SC(\Sigma_1): \sqrt n+\sqrt m\not=\sqrt k+\sqrt l,
1\leq N\leq M\leq y, 1\leq L\leq K\leq y,\\
&&\hspace{20mm}N\leq L, n\sim N, m\sim M, k\sim K, l\sim L.
\end{eqnarray*}

If $M\geq 100K, $ then $|\sqrt n+\sqrt m-\sqrt k-\sqrt l|\gg
M^{1/2},$ so the trivial estimate yields
$$G(N,M,K,L) \ll \frac{T^{\frac{3}{2}+\varepsilon}NMKL}{(NMKL)^{3/4}M^{1/2}}
\ll T^{\frac{3}{2}+\varepsilon} y^{\frac{1}{2}}\ll
T^{15/8+\varepsilon}.$$

If $K>100M,$ we get the same estimate. So later we always suppose
that $M\asymp K.$

Let $\eta=\sqrt n+\sqrt m-\sqrt k-\sqrt l.$ Write
\begin{equation}
G(N,M,K,L,R)=G_1+G_2+G_3,
\end{equation}
where
\begin{eqnarray*}
&&G_1:=T^2\sum_{|\eta|\leq T^{-1/2}}g,\\
&&G_2:=T^{\frac{3}{2}}\sum_{ T^{-1/2}<|\eta|\leq 1}g|\eta|^{-1},\\
&&G_3:=T^{\frac{3}{2}}\sum_{ |\eta|\gg 1}g|\eta|^{-1}.
\end{eqnarray*}

We estimate $G_1$ first.
 From $|\eta|\leq T^{-1/2}$ we get $M\asymp K\gg T^{1/7}$ via Lemma 2.
By Lemma 5 we get
\begin{align}
G_1&\ll \frac{T^{2+\varepsilon}}{(NMKL)^{3/4}} {\cal
A}_1(N,M,K,L;T^{-1/2})\\
&\ll \frac{T^{2+\varepsilon}}{(NMKL)^{3/4}}
\left(T^{-1/2}K^{1/2}NML+NLK^{1/2}  \right)\nonumber\\
&\ll
T^{3/2+\varepsilon}(NL)^{1/4}+T^{2+\varepsilon}(NL)^{1/4}K^{-1}\nonumber\\
&\ll
T^{3/2+\varepsilon}y^{1/2}+T^{2+\varepsilon}(NL)^{1/4}K^{-1}\nonumber\\
&\ll
T^{15/8+\varepsilon}+T^{2+\varepsilon}(NL)^{1/4}K^{-1}.\nonumber
\end{align}
By Lemma 3 we get(notice $N\leq L\leq K$)
\begin{align}
G_1&\ll \frac{T^{2+\varepsilon}}{(NMKL)^{3/4}} {\cal
A}_{-}(N,M,K,L;T^{-1/2})\\
&\ll \frac{T^{2+\varepsilon}}{(NMKL)^{3/4}}
\left(T^{-1/8}N^{7/8}+N^{1/2}\right)\left(T^{-1/8}L^{7/8}+L^{1/2}\right)
\left(T^{-1/4}K^{7/4}+K\right)\nonumber\\
&\ll
T^{2+\varepsilon}\left(T^{-1/8}N^{1/8}+N^{-1/4}\right)\left(T^{-1/8}L^{1/8}+L^{-1/4}\right)
\left(T^{-1/4}K^{1/4}+K^{-1/2}\right)\nonumber\\
&\ll
T^{2+\varepsilon}\left(T^{-1/4}(NL)^{1/8}+T^{-1/8}L^{1/8}N^{-1/4}+(NL)^{-1/4}\right)
\left(T^{-1/4}K^{1/4}+K^{-1/2}\right)\nonumber\\
&\ll
T^{2+\varepsilon}T^{-1/4}(NL)^{1/8}\left(T^{-1/4}K^{1/4}+K^{-1/2}\right)\nonumber\\
&\hspace{3mm}+T^{2+\varepsilon}\left(T^{-1/8}L^{3/8}(NL)^{-1/4}+(NL)^{-1/4}\right)
\left(T^{-1/4}K^{1/4}+K^{-1/2}\right)\nonumber\\
&\ll
T^{3/2+\varepsilon}y^{1/2}+T^{7/4+\varepsilon}K^{-1/4}\nonumber\\&\hspace{3mm}
+ T^{2+\varepsilon}\left(T^{-1/4}K^{1/4}+K^{-1/2}\right)
\left(T^{-1/8}K^{3/8}+1\right)(NL)^{-1/4}\nonumber\\
&\ll
T^{15/8+\varepsilon}+T^{2+\varepsilon}K^{-1/2}\left(T^{-1/4}K^{3/4}+1\right)
\left(T^{-1/8}K^{3/8}+1\right)(NL)^{-1/4}\nonumber\\
&\ll T^{15/8+\varepsilon}+T^{2+\varepsilon}K^{-1/2}
\left(T^{-3/8}K^{9/8}+1\right)(NL)^{-1/4}\nonumber.\nonumber
\end{align}

From (3.8) and (3.9) we get
\begin{align}
G_1&\ll
T^{15/8+\varepsilon}+T^{2+\varepsilon}\min\left((NL)^{1/4}K^{-1},K^{-1/2}
\left(T^{-3/8}K^{9/8}+1\right)(NL)^{-1/4}\right)\\
&\ll
T^{15/8+\varepsilon}+T^{2+\varepsilon}\left((NL)^{1/4}K^{-1}\right)^{1/2}\left(K^{-1/2}
\left(T^{-3/8}K^{9/8}+1\right)(NL)^{-1/4}\right)^{1/2}\nonumber\\
&\ll
T^{15/8+\varepsilon}+T^{2+\varepsilon}K^{-3/4}(T^{-3/16}K^{9/16}+1)\nonumber\\
&\ll T^{15/8+\varepsilon}+T^{2+\varepsilon}K^{-3/4}\ll
T^{53/28+\varepsilon}\nonumber
\end{align}
if we notice $K\gg T^{1/7}.$

 Now we estimate $G_2.$  By a splitting argument
we get that the estimate
\begin{equation}
G_2\ll
\frac{T^{3/2+\varepsilon}}{(NMKL)^{3/4}\delta}\sum_{\stackrel{\delta<|\eta|\leq
2\delta }{\eta\not= 0}}1
\end{equation}
holds for some $T^{-1/2}\leq \delta\leq 1.$ By Lemma 5  we get that
\begin{align}
G_2&\ll \frac{T^{3/2+\varepsilon}}{(NMKL)^{3/4}\delta}{\cal
A}_1(N,M,K,L;2\delta)\\
&\ll \frac{T^{3/2+\varepsilon}}{(NMKL)^{3/4}\delta}
(\delta K^{1/2}NML+NLK^{1/2})\nonumber\\
&\ll  T^{3/2+\varepsilon}y^{1/2}+
T^{3/2+\varepsilon}(K\delta)^{-1}(NL)^{1/4} \nonumber\\
&\ll
T^{15/8+\varepsilon}+T^{3/2+\varepsilon}(K\delta)^{-1}(NL)^{1/4}\nonumber.
\end{align}
By Lemma 3 we get(notice $N\leq L\leq K$)
\begin{align}
G_2&\ll \frac{T^{3/2+\varepsilon}}{(NMKL)^{3/4}\delta}{\cal
A}_{-}(N,M,K,L;2\delta)\\
&\ll \frac{T^{3/2+\varepsilon}}{(NMKL)^{3/4}\delta}
(\delta^{1/4}N^{7/8}+N^{1/2})(\delta^{1/4}L^{7/8}+L^{1/2})
(\delta^{1/2}K^{7/4}+K)\nonumber\\
&\ll
T^{3/2+\varepsilon}(N^{1/8}+N^{-1/4}\delta^{-1/4})(L^{1/8}+L^{-1/4}\delta^{-1/4})
(K^{1/4}+K^{-1/2}\delta^{-1/2})\nonumber\\
&\ll
T^{3/2+\varepsilon}\left((NL)^{1/8}+L^{1/8}N^{-1/4}\delta^{-1/4}+(NL)^{-1/4}\delta^{-1/2}\right)
(K^{1/4}+K^{-1/2}\delta^{-1/2})\nonumber\\
&\ll
T^{3/2+\varepsilon}(NL)^{1/8}K^{1/4}+T^{3/2+\varepsilon}(NL)^{1/8}K^{-1/2}\delta^{-1/2}\nonumber\\
&\hspace{3mm}+T^{3/2+\varepsilon}(K^{1/4}+K^{-1/2}\delta^{-1/2})(L^{3/8}\delta^{1/4}+1)(NL)^{-1/4}\delta^{-1/2}\nonumber\\
&\ll T^{3/2+\varepsilon}y^{1/2}+T^{3/2+\varepsilon}\delta^{-1/2} +
T^{3/2+\varepsilon}K^{-1/2}\delta^{-1}(K^{3/4}\delta^{1/2}+1)(K^{3/8}\delta^{1/4}+1)(NL)^{-1/4}\nonumber\\
&\ll T^{15/8+\varepsilon}+
T^{3/2+\varepsilon}K^{-1/2}\delta^{-1}(K^{9/8}\delta^{3/4}+1)(NL)^{-1/4},\nonumber
\end{align}
 where the bound  $\delta\gg T^{-1/2}$ was used  to the term
$T^{3/2+\varepsilon}\delta^{-1/2}.$

From (3.12) and (3.13) we get
\begin{align}
G_2&\ll
T^{15/8+\varepsilon}+\frac{T^{3/2+\varepsilon}}{\delta}\min\left(\frac{(NL)^{1/4}}{K},
\frac{K^{9/8}\delta^{3/4}+1}{K^{1/2}(NL)^{1/4}}\right)\\
&\ll T^{15/8+\varepsilon}+\frac{T^{3/2+\varepsilon}}{\delta}
\left(\frac{(NL)^{1/4}}{K}\right)^{1/2}
\left(\frac{K^{9/8}\delta^{3/4}+1}{K^{1/2}(NL)^{1/4}}\right)^{1/2}\nonumber\\
&\ll
T^{15/8+\varepsilon}+T^{3/2+\varepsilon}\delta^{-1}K^{-3/4}(K^{9/16}\delta^{3/8}+1).\nonumber
\end{align}

If $\delta\gg K^{-3/2},$ then (3.14) implies(recall $\delta\gg
T^{-1/2}$)
\begin{equation}
G_2\ll
T^{15/8+\varepsilon}+T^{3/2+\varepsilon}K^{-3/16}\delta^{-5/8}\ll
T^{15/8+\varepsilon}.
\end{equation}
If $\delta\ll K^{-3/2},$ then (3.14) becomes
\begin{equation}
G_2\ll T^{15/8+\varepsilon}+T^{3/2+\varepsilon}\delta^{-1}K^{-3/4}.
\end{equation}
Since we have $\delta\gg K^{-7/2}$ by Lemma 2 and $\delta\gg
T^{-1/2},$ we get $\delta^{-1}\ll \min(K^{7/2},T^{1/2})$ and thus
from (3.16) we get
\begin{align}
G_2&\ll
T^{15/8+\varepsilon}+\min(T^{2+\varepsilon}K^{-3/4},T^{3/2+\varepsilon}K^{11/4})\\
&\ll
T^{15/8+\varepsilon}+(T^{2+\varepsilon}K^{-3/4})^{11/14}(T^{3/2+\varepsilon}K^{11/4})^{3/14}\nonumber\\
&\ll T^{53/28+\varepsilon}.\nonumber
\end{align}

For $G_3,$ by a splitting argument and Lemma 5 again (notice
$|\eta|\gg 1 $) we get
\begin{equation}
G_3\ll \frac{T^{3/2+\varepsilon}}{(NMKL)^{3/4}\delta}
\sum_{\delta<|\eta|\leq 2\delta,\delta\gg 1}1
\end{equation}
$$\ll \frac{T^{3/2+\varepsilon}}{(NMKL)^{3/4}}K^{1/2}NML
\ll  T^{\frac{3}{2}+\varepsilon}y^{\frac{1}{2}}\ll
T^{15/8+\varepsilon}.$$

Combining (3.6), (3.7), (3.10) and (3.15)-(3.18) we get
\begin{equation}
\int_T^{2T}S_2(x)dx\ll T^{53/28+\varepsilon}.
\end{equation}

In the same way we can show that
\begin{equation}
\int_T^{2T}S_3(x)dx\ll T^{53/28+\varepsilon}
\end{equation}
by  Lemma 3 and  Lemma 6.

From (3.2)-(3.5) , (3.19) and (3.20) we get
\begin{equation}
\int_T^{2T}\Delta^4(x)dx=\frac{3c_2}{32\pi^4}\int_T^{2T}xdx+O(T^{53/28+\varepsilon}),
\end{equation}
which implies Theorem 1 immediately.

\section{\bf Preliminary Lemmas for Theorem 2}

In order to prove Theorem 2 , we need the following Lemmas.

{\bf Lemma 7.} We have
$$E(t)=\Sigma_1(t)+\Sigma_2(t)+O(\log^2 t)$$
with
\begin{eqnarray}
&&\Sigma_1(t):=\frac{1}{\sqrt 2}\sum_{n\leq N}h(t,n)\cos(f(t,n)),\\
&&\Sigma_2(t):=-2\sum_{n\leq
N^{\prime}}d(n)n^{-1/2}(\log\frac{t}{2\pi n})^{-1}
\cos(t\log\frac{t}{2\pi n}-t+\frac{\pi}{4}),\\
&&h(t,n):=(-1)^nd(n)n^{-1/2}(\frac{t}{2\pi n}+\frac{1}{4})^{-1/4}
(g(t,n))^{-1},\\
&&g(t,n):=arsinh((\frac{\pi n}{2t})^{1/2}),\\
&&f(t,n):=2tg(t,n)+(2\pi nt+\pi^2n^2)^{1/2}-\pi/4,\\
&&\quad At\leq N\leq A^{\prime}t, N^{\prime}:=
t/2\pi+N/2-(N^2/4+Nt/2\pi)^{1/2},
\end{eqnarray}
where $0<A<A^{\prime}$ are any fixed constants.

\begin{proof}This is the famous Atkinson's formula, see Atkinson [1] or Ivi\'c[5, Theorem 15.1].
\end{proof}

{\bf Lemma 8.} Suppose $Y>1.$ Define
\begin{eqnarray*}
c_2^{*}:&&=\sum_{\sqrt{n}+\sqrt{m}=\sqrt{k}+\sqrt{l}}\frac{(-1)^{n+m+k+l}
d(n)d(m)d(k)d(l)}{(nmkl)^{3/4}},\\
c_2^{*}(Y):&&=\sum_{\stackrel{\sqrt{n}+\sqrt{m}=\sqrt{k}+\sqrt{l}}{n,m,k,l\leq
Y}} \frac{(-1)^{n+m+k+l}d(n)d(m)d(k)d(l)}{(nmkl)^{3/4}},\\
c_2(Y):&&=\sum_{\stackrel{\sqrt{n}+\sqrt{m}=\sqrt{k}+\sqrt{l}}{n,m,k,l\leq
Y}} \frac{d(n)d(m)d(k)d(l)}{(nmkl)^{3/4}}.
\end{eqnarray*}
Then we have
\begin{eqnarray*}
c_2=c_2^{*},\quad  c_2(Y)=c_2^{*}(Y),\quad  |c_2-c_2(Y)|\ll
Y^{-1/2+\varepsilon}.
\end{eqnarray*}

\begin{proof}
The estimate $|c_2-c_2(Y)|\ll Y^{-1/2+\varepsilon}$ is a special
case of Lemma 3.1 of  the author[14] . The equalities
 $c_2=c_2^{*}$ and $c_2(Y)=c_2^{*}(Y)$
follow from the fact that if
$\sqrt{n_1}+\sqrt{n_2}=\sqrt{n_3}+\sqrt{n_4},$ then
$n_1+n_2+n_3+n_4$ must be an even number .

\end{proof}

{\bf Lemma 9.} Suppose $Y>1,$ then we have
$$H_1(Y):=\sum_{\stackrel{\sqrt{n}+\sqrt{m}=\sqrt{k}+\sqrt{l}}{n,m,k,l\leq Y}}
\frac{d(n)d(m)d(k)d(l)\max(n,m,k,l)^{3}}{(nmkl)^{3/4}}\ll
Y^{5/2+\varepsilon}.$$

\begin{proof}
Suppose $\sqrt{n}+\sqrt{m}=\sqrt{k}+\sqrt{l},$ then we have

 (1)  $n=k, m=l \hspace{2mm}\mbox{or}\hspace{2mm} n=l, m=k;
 $

 or

$(2) n\not= k,l.$

If the case (2) holds, then by  a classical result of Besicotitch,
we know that
$$n=n_1^2h,m=m_1^2h,k=k_1^2h,l=l_1^2h, n_1+m_1=k_1+l_1,
\mu(h)\not= 0.$$ Thus we get
\begin{eqnarray*}
H_1(Y)&&\ll \Sigma_1+\Sigma_2,\\
 \Sigma_1&&\ll \sum_{n,k\leq
Y}\frac{d^2(n)d^2(m)\max(n,k)^{3}}{(nk)^{3/2}}\ll Y^{5/2}\log^3 Y,\\
\Sigma_2&&\ll Y^{\varepsilon}\sum_{h<
Y}\sum_{\stackrel{n_1+m_1=k_1+l_1}{n_1,m_1,k_1,l_1\leq
Y^{1/2}h^{-1/2}}}\frac{\max(n_1,m_1,k_1,l_1)^6}{(n_1m_1k_1l_1)^{3/2}}\\
&&\ll Y^{\varepsilon}\sum_{h<
Y}\sum_{\stackrel{n_1+m_1=k_1+l_1}{n_1,m_1,l_1\leq k_1\leq
Y^{1/2}h^{-1/2}}}\frac{k_1^{9/2}}{(n_1m_1l_1)^{3/2}}\\
&&\ll Y^{\varepsilon}\sum_{h<
Y}\sum_{l_1}l_1^{-3/2}\sum_{\stackrel{n_1+m_1>k_1}{n_1,m_1\leq
k_1\leq Y^{1/2}h^{-1/2}}}\frac{k_1^{9/2}}{(n_1m_1)^{3/2}}\\
&&\ll Y^{\varepsilon}\sum_{h<
Y}\sum_{l_1}l_1^{-3/2}\sum_{n_1}n_1^{-3/2}\sum_{k_1\ll m_1\leq k_1\leq Y^{1/2}h^{-1/2}}k_1^3\\
&&\ll Y^{\varepsilon}\sum_{h<Y}(Y^{1/2}h^{-1/2})^5\ll
Y^{5/2+\varepsilon}.
\end{eqnarray*}
\end{proof}

{\bf Lemma 10.} Suppose $Y>1,$ then we have
$$H_2(Y):=\sum_{\stackrel{\sqrt{n}+\sqrt{m}+\sqrt{k}=\sqrt{l}}{n, m, k, l\leq Y}}
\frac{d(n)d(m)d(k)d(l)l^{3/4}}{(nmk)^{3/4}}\ll
Y^{1/2+\varepsilon}.$$

\begin{proof}
If $\sqrt{n}+\sqrt{m}+\sqrt{k}=\sqrt{l},$ then we have
$$n=n_1^2h,m=m_1^2h,k=k_1^2h,l=l_1^2h, n_1+m_1+k_1=l_1,
\mu(h)\not= 0.$$ Thus we get
\begin{eqnarray*}
H_2(Y) &&\ll Y^{\varepsilon}\sum_{h(n_1+m_1+k_1)^2\leq
Y}\frac{(n_1+m_1+k_1)^{3/2}}{h^{3/2}(n_1m_1k_1)^{3/2}}\\
&&\ll Y^{\varepsilon}\sum_{h<Y}h^{-3/2}\sum_{n_1\leq m_1\leq k_1\leq
Y^{1/2}h^{-1/2}}n_1^{-3/2}m_1^{-3/2}\ll Y^{1/2+\varepsilon}.
\end{eqnarray*}
\end{proof}

{\bf Lemma 11.} Suppose  $f_j(t)(1\leq j\leq k)$ and $g(t)$ are
 continuous , monotonic real-valued functions on $[a,b]$
 and let $g(t)$ have a continuous , monotonic derivative on $[a,b].$
If $| f_j(t)|\leq A_j (1\leq j\leq k),
 |g^{\prime}(t)|\gg \Delta$ for any $t\in [a,b],$ then
$$\int_a^bf_1(t)\cdots f_k(t)e(g(t))dt\ll A_1\cdots A_k\Delta^{-1}.$$

\begin{proof}This is Lemma 15.3 of Ivi\'c[5].
\end{proof}

\section{\bf Proof of Theorem 2}

Suppose $T\geq 10.$ It suffices to evaluate $\int_T^{2T}E^4(t)dt.$
Let  $y:=T^{1/3-\varepsilon}$ . For any $T\leq t\leq 2T,$ define
\begin{eqnarray*}
&&{\cal E}_1(t):=\frac{1}{\sqrt 2}\sum_{n\leq y}h(t,n)\cos(f(t,n)),
\quad {\cal E}_2(t):=E(t)-{\cal E}_1(t).
\end{eqnarray*}

From  the inequality $(a+b)^4-a^4\ll |b|^3|a|+|b|^4,$ we get
\begin{equation}
\int_T^{2T}E^4(t)dt=\int_T^{2T}{\cal E}_1^4(t)dt
+O(\int_T^{2T}|{\cal E}_1(t)|^3|{\cal E}_2(t)|dt)
+O(\int_T^{2T}|{\cal E}_2(t)|^4dt).
\end{equation}

\subsection{\bf Evaluation of $\int_T^{2T}{\cal E}_1^4(t)dt$}\

In this subsection , we shall evaluate the integral
$\int_T^{2T}{\cal E}_1^4(t)dt.$ Similar to Tsang[11], we can write
\begin{equation}
{\cal
E}_1^4(t)=\frac{3}{32}S_5(t)+\frac{3}{32}S_6(t)+\frac{1}{8}S_7(t)
+\frac{1}{8}S_8(t)+\frac{1}{32}S_9(t),
\end{equation}
where
\begin{eqnarray*}
&&S_5(t):=\sum_{\stackrel{n,m,k,l\leq y}{\sqrt n+\sqrt m=\sqrt
k+\sqrt
l}}H(t;n,m,k,l)\cos(F_1(t;n,m,k,l)),\\
&&S_6(t):=\sum_{\stackrel{n,m,k,l\leq y}{\sqrt n+\sqrt m\not=\sqrt
k+\sqrt
l}}H(t;n,m,k,l)\cos(F_1(t;n,m,k,l)),\\
&&S_7(t):=\sum_{\stackrel{n,m,k,l\leq y}{\sqrt n+\sqrt m+\sqrt
k=\sqrt
l}}H(t;n,m,k,l)\cos(F_2(t;n,m,k,l)),\\
&&S_8(t):=\sum_{\stackrel{n,m,k,l\leq y}{\sqrt n+\sqrt m+\sqrt
k\not=\sqrt
l}}H(t;n,m,k,l)\cos(F_2(t;n,m,k,l)),\\
&&S_9(t):=\sum_{n,m,k,l\leq
y}H(t;n,m,k,l)\cos(F_3(t;n,m,k,l)),\\
&&H(t;n,m,k,l):=h(t,n)h(t,m)h(t,k)h(t,l),\\
&&F_1(t;n,m,k,l):=f(t,n)+f(t,m)-f(t,k)-f(t,l),\\
&&F_2(t;n,m,k,l):=f(t,n)+f(t,m)+f(t,k)-f(t,l),\\
&&F_3(t;n,m,k,l):=f(t,n)+f(t,m)+f(t,k)+f(t,l).
\end{eqnarray*}

We first estimate  the integral $\int_T^{2T}S_5(t)dt.$ For $n\leq
y,$ it is  easy to check that
\begin{eqnarray}
&&h(t,n)=\frac{2^{3/4}}{\pi^{1/4}}\frac{(-1)^nd(n)}{n^{3/4}}t^{1/4}
(1+O(\frac{n}{t})),\\
&& f(t,n)=2^{3/2}(\pi nt)^{1/2}-\pi/4+O(n^{3/2}t^{-1/2}),\\
&&f^{\prime}(t,n)=2^{1/2}(\pi n)^{1/2}t^{-1/2}+O(n^{3/2}t^{-3/2}).
\end{eqnarray}
If $ \sqrt n+\sqrt m=\sqrt k+\sqrt l$, then
\begin{equation}
\cos(F_1(n,m,k,l))=\cos(O(\frac{D^{3/2}}{t^{1/2}}))=
1+O(\frac{D^{3}}{t}),
\end{equation}
where $D:=\max(n,m,k,l).$
 So from (5.3), (5.6), Lemma 8 and Lemma 9 we get
\begin{align}
\int_T^{2T}S_5(t)dt&=\sum_{\stackrel{n,m,k,l\leq y}{\sqrt n+\sqrt
m=\sqrt k+\sqrt l}}\int_T^{2T}H(t;n,m,k,l)\cos(F_1(t;n,m,k,l))dt\\
&=\frac{8}{\pi}\sum_{\stackrel{n,m,k,l\leq y}{\sqrt n+\sqrt m=\sqrt
k+\sqrt
l}}\frac{(-1)^{n+m+k+l}d(n)d(m)d(k)d(l)}{(nmkl)^{3/4}}\nonumber\\
&\hspace{3mm}\times
\int_T^{2T}t(1+O(\frac{D}{t}))(1+(\frac{D^{3}}{t}))dt\nonumber\\
&=\frac{8}{\pi}\sum_{\stackrel{n,m,k,l\leq y}{\sqrt n+\sqrt m=\sqrt
k+\sqrt
l}}\frac{(-1)^{n+m+k+l}d(n)d(m)d(k)d(l)}{(nmkl)^{3/4}}\nonumber\\
&\hspace{3mm}\times
\int_T^{2T}t(1+(\frac{D^{3}}{t}))dt\nonumber\\
&=\frac{8}{\pi}\sum_{\stackrel{n,m,k,l\leq y}{\sqrt n+\sqrt m=\sqrt
k+\sqrt
l}}\frac{(-1)^{n+m+k+l}d(n)d(m)d(k)d(l)}{(nmkl)^{3/4}}\int_T^{2T}tdt\nonumber\\
&\hspace{10mm}+O(TH_1(y))\nonumber\\
&=\frac{8c_2}{\pi}\int_T^{2T}tdt+O(T^{1+\varepsilon}y^{5/2}+T^{2+\varepsilon}y^{-1/2})\nonumber\\
&=\frac{8c_2}{\pi}\int_T^{2T}tdt+O(T^{11/6+\varepsilon}).\nonumber
\end{align}

Now we estimate $\int_T^{2T}S_6(t)dt.$ From (5.5) we get
$$F_1^{\prime}(t;n,m,k,l)=(2\pi)^{1/2}\eta t^{-1/2}+O(D^{3/2}t^{-3/2}), $$
where $\eta=\sqrt n+\sqrt m-\sqrt k-\sqrt l.$ Write
\begin{equation}
\int_T^{2T}S_6(t)dt=\int_{|\eta|\leq
T^{-1/2}}S_6(t)dt+\int_{|\eta|>T^{-1/2}}S_6(t)dt.
\end{equation}
 If $|\eta|\leq T^{-1/2},$ the by (5.3) and the trivial
estimate we get
\begin{equation}
\int_{|\eta|\leq T^{-1/2}}S_6(t)dt\ll T^2\sum_{\stackrel{n,m,k,l\leq
y;|\eta|\leq T^{-1/2}}{\sqrt n+\sqrt m\not=\sqrt k+\sqrt
l}}\frac{d(n)d(m)d(k)d(l)}{(nmkl)^{3/4}}.
\end{equation}

If $|\eta|> T^{-1/2},$ then $|F_1^{\prime}(t;n,m,k,l)|\gg
|\eta|T^{-1/2},$ thus from (5.3) and Lemma 11 we get
\begin{equation}
\int_{|\eta|> T^{-1/2}}S_6(t)dt\ll
T^{3/2}\sum_{\stackrel{n,m,k,l\leq y;|\eta|> T^{-1/2}}{\sqrt n+\sqrt
m\not=\sqrt k+\sqrt l}}\frac{d(n)d(m)d(k)d(l)}{(nmkl)^{3/4}|\eta|}.
\end{equation}

From (5.9), (5.10) and the estimate in Section 3 we get
\begin{equation}
\int_T^{2T}S_6(t)dt\ll \sum_{\stackrel{n,m,k,l\leq y}{\sqrt n+\sqrt
m\not=\sqrt k+\sqrt
l}}\frac{d(n)d(m)d(k)d(l)}{(nmkl)^{3/4}}\min(T^2,T^{3/2}|\eta|^{-1})\ll
T^{53/28+\varepsilon}.
\end{equation}

 If $\sqrt n+\sqrt m+\sqrt k=\sqrt l,$ then from (5.4) we have
\begin{eqnarray*}
&&F_2(t;n,m,k,l)=-\pi/2+O(l^{3/2}t^{-1/2}),
\hspace{2mm}\cos{F_2(t;n,m,k,l)}\ll l^{3/2}t^{-1/2}.
\end{eqnarray*}
Thus from (5.3), the trivial estimate and Lemma 10 we get
\begin{equation}
\int_T^{2T}S_7(t)dt\ll T^{3/2}H_2(y)\ll
T^{3/2}y^{1/2+\varepsilon}\ll T^{5/3+\varepsilon}.
\end{equation}

Similar to the integral $\int_T^{2T}S_6(t)dt,$
 we have
 \begin{equation}
\int_T^{2T}S_8(t)dt\ll T^{53/28+\varepsilon}.
 \end{equation}

From (5.5) we get
$$F_3^{\prime}(t;n,m,k,l)\gg (\sqrt n+\sqrt m+\sqrt k+\sqrt l)T^{-1/2},$$
which combining  (5.3) and Lemma 11 implies
\begin{equation}
\int_T^{2T}S_9(t)dt\ll
\sum_{n,m,k,l}\frac{d(n)d(m)d(k)d(l)T^{3/2}}{(nmkl)^{3/4}(\sqrt
n+\sqrt m+\sqrt k+\sqrt l)}\ll T^{3/2+\varepsilon}y^{1/2}\ll
T^{5/3+\varepsilon}.
\end{equation}

From (5.2), (5.7),  (5.11)-(5.14) we get
\begin{equation}
\int_T^{2T}{\cal E}_1^4(t)dt=\frac{3c_2}{4\pi}\int_T^{2T}tdt
+O(T^{53/28+\varepsilon}).
\end{equation}

\subsection{\bf Completion of proof of Theorem 2}\

Let $A_0=35/8.$ Ivi\'c[5, Thm 15.7] proved the estimate
\begin{equation}
\int_1^T|E(t)|^{A_0}dt\ll T^{1+A_0/4+\varepsilon} .
\end{equation}
By his method we can show
\begin{equation}
\int_T^{2T}|{\cal E}_1(t)|^{A_0}dt\ll T^{1+A_0/4+\varepsilon} .
\end{equation}
Thus
\begin{equation}
\int_T^{2T}|{\cal E}_2(t)|^{A_0}dt\ll T^{1+A_0/4+\varepsilon} .
\end{equation}
We also have
\begin{equation}
\int_T^{2T}|{\cal E}_2(t)|^2dt\ll T^{3/2+\varepsilon}y^{-1/2},
\end{equation}
which is the formula (4.15) of the author[14]. From (5.18) , (5.19)
and the H\"older's inequality we get  that the estimate
\begin{equation}
\int_T^{2T}|{\cal E}_2(t)|^Adt\ll
T^{1+A/4+\varepsilon}y^{-(A_0-A)/2(A_0-2)}
\end{equation}
holds for  any $2<A<A_0.$ The details of the above estimates can be
found in the author[14].

From (5.17), (5.20) and the H\"older's inequality we get
\begin{align}
\int_T^{2T}|{\cal E}_1^3(t){\cal E}_2(t)|dt&\ll \left(\int_1^T|{\cal
E}_1(t)|^{A_0}dt\right)^{3/A_0}\left(\int_1^T|{\cal
E}_2(t)|^{A_0/(A_0-3)}dt\right)^{(A_0-3)/A_0}\\
&\ll T^{2+\varepsilon}y^{-(A_0-4)/2(A_0-2)}\ll
T^{2-19/108+\varepsilon}\nonumber.
\end{align}

From (5.1), (5.15), (5.20) with $A=4$ and (5.21)  we get
\begin{equation}
\int_T^{2T}E^4(t)dt=\frac{3c_2}{4\pi}\int_T^{2T}tdt
+O(T^{53/28+\varepsilon})
\end{equation}
and Theorem 2 follows.


\begin{thebibliography}{99}
 \bibitem{s1}F. V. Atkinson, the mean value of the Riemann zeta-function,  Acta
Math. 1949, 81: 353-376.

\bibitem{s2}Cai Yingchun, On the third and fourth power moments of Fourier
coefficients of cusp forms.  Acta Math. Sinica (N.S.) 13 (1997), no.
4, 443--452 .

\bibitem{s2}M. N. Huxley,  Exponential sums and
Lattice points III, Proc. London Math. Soc., Vol.{\bf 87}(3)(2003),
591-609.

\bibitem{s3}A. Ivi\'c, On some problems involving the mean square of
 $\zeta(\frac 12+it)$. Bull. Cl. Sci. Math. Nat. Sci. Math. 1998, 23: 71--76.

\bibitem{s4}A. Ivi\'c, The Riemann zeta-function. John. Wiley and Sons, 1985.

\bibitem{s5}A. Ivi\'c, Lectures on mean values of the Riemann zeta-function,
Lectures On Math. and Physics 82, Tata Inst. Fund. Res., Bombay,
1991.
\bibitem{IS}A. Ivi\'c  and P. Sargos, On the higher power moments
of the error term in the divisor problem, Illinois Journal of Math.
{\bf 51}(2)(2007), 353-377.


\bibitem{s7}M. Jutila, On a formula of Atkinson.
Topics in classical number theory, Vol. I, II (Budapest, 1981),
807--823, Colloq. Math. Soc. Jnos Bolyai, 34, North-Holland,
Amsterdam,
 1984.

 \bibitem{s2}O. Robert and P. Sargos, Three-dimemsional
 exponential sums with monomials, J. Reine angew. Math. 591(2006), 1-20.



\bibitem{s8}K. C. Tong, On divisor problem III, Acta math. Sinica {\bf 6}
(1956), 515-541.

\bibitem{s9}Kai-Man Tsang, Higher-power moments of $\Delta(x), E(t)$ and $P(x)$,
Proc. London Math. Soc.(3){\bf 65}(1992), 65-84.

\bibitem{s10}Wenguang Zhai , On higher-power moments of $\Delta(x),$ Acta
Arith. Vol.{\bf 112}(2004), 1-24.


\bibitem{s11} Wenguang Zhai, On higher-power moments of $\Delta(x) (II),$ Acta Arith.
Vol.{\bf 114}(2004), 35-54.

\bibitem{s12} Wenguang Zhai, On higher-power moments of $E(t),$ Acta
Arith. Vol.{\bf 115}(2004),  329-348.
\end{thebibliography}
\end{document}